\documentclass[11pt]{article}
\usepackage{a4wide}

\usepackage[a4paper, margin=2.3cm]{geometry}
\usepackage{amsmath,amssymb,amsthm}
\usepackage{color}
\usepackage{parskip}
\usepackage{hyperref}
\usepackage{parskip}
\usepackage{setspace}
\usepackage{graphicx}

\newtheorem{theorem}{Theorem}[section]
\newtheorem{lemma}[theorem]{Lemma}

\newtheorem{corollary}[theorem]{Corollary}

\newcommand{\x} {\pmb x}
\newcommand{\y} {\pmb y}
\newcommand{\R}{\mathrm{\mathcal R}}
\newcommand{\A}{\mathrm{\mathcal A}}
\newcommand{\B}{\mathrm{\mathcal B}}
\newcommand{\C}{\mathrm{\mathcal C}}
\newcommand{\D}{\mathrm{\mathcal D}}
\newcommand{\E}{\mathrm{\mathcal E}}

\title{A note on sum-product estimates over finite valuation rings}
\author{Pham Duc Hiep}

\begin{document}
\date{}
\maketitle

\begin{abstract}
Let $\R$ be a finite valuation ring of order $q^r$ with $q$ a power of an odd prime number, and $\A$ be a set in $\R$. 
In this paper, we improve a recent result due to Yazici (2018) on a sum-product type problem. More precisely, we will prove that 
\bigskip
\begin{enumerate}
\item If $|\A|\gg q^{r-\frac{1}{3}}$, then 
\[\max\left\lbrace |\A+\A|, |\A^2+\A^2|\right\rbrace \gg  q^{\frac{r}{2}}|\A|^{\frac{1}{2}}.\]
\item If $q^{r-\frac{3}{8}}\ll |\A|\ll q^{r-\frac{1}{3}}$, then
\[\max\left\lbrace |\A+\A|, |\A^2+\A^2|\right\rbrace \gg  \frac{|\A|^2}{q^{\frac{2r-1}{2}}}.\]
\item If $|\A+\A||\A|^2\gg q^{3r-1}$ and $2q^{r-1}\le |\A|\ll q^{r-\frac{3}{8}}$, then
\[\max\left\lbrace |\A+\A|, |\A^2+\A^2|\right\rbrace \gg  q^{r/3}|\A|^{2/3}.\]
\end{enumerate}
\end{abstract}
\def\thefootnote{\empty}

\section{Introduction}\label{section1}
Let $A$ be a subset of integers, we define the sum set and the product set, respectively, as follows 
$$\A+\A = \{a+b\colon a, b \in \A\},$$
$$\A\cdot \A=\{a\cdot b\colon  a,b \in \A\}.$$

In 1983, Erd\H{o}s and Szemer\'edi \cite{es} proved that there is no set $\A$ which is both highly additive structured and multiplicative structured at the same time.  More precisely, they proved that
$$\max\{ |\A+\A|, |\A \cdot \A|\} \geq |\A|^{1+c},$$
for some constant $c >0$

Over last twenty years, there was intensive progress on improving the constant $c$ and on studying variants in different settings, for example, over finite fields and the complex numbers. We refer the interested reader to \cite{mot, garaev, rx, RSS, shen, vu} and references therein for more details. 

Let $\R$ be a finite valuation ring, i.e. a finite, local and principal ring.
The first result on sum-product type problems in the context of finite valuation rings was given in \cite{ham} by Ham, Pham and Vinh. In particular, they proved the following theorems. 
\bigskip

\begin{theorem}[Ham-Pham-Vinh, \cite{ham}]
Let $\R$ be a finite valuation ring of order $q^r$, and $\R^*$ denote the set of units of $\R$. Let $G$ be a subgroup of $\R^*$, and $f(x,y)=g(x)(h(x)+y)$ be defined on $G\times \R^*$, where $g,h:G\rightarrow\R^*$ are arbitrary functions. Put $m=\displaystyle\max_{t\in \R}|\{x\in G: g(x)h(x)=t\}|$. For any sets $\A\subset G$ and $\B,\C\subset \R^*$, we have 
$$|f(\A,\B)||\B\cdot \C|\gg \min\Bigg\{\dfrac{q^r|B|}m,\dfrac{|\A||\B|^2|\C|}{m^2q^{2r-1}}\Bigg\}.$$
\end{theorem}
\begin{theorem}[Ham-Pham-Vinh, \cite{ham}]\label{ham}
Let $\R$ be a finite valuation ring of order $q^r$, and $\R^*$ denote the set of units of $\R$. Let $G$ be a subgroup of $\R^*$, and $f(x,y)=g(x)(h(x)+y)$ be defined on $G\times \R^*$, where $g,h:G\rightarrow\R^*$ are arbitrary functions. Put $m=\displaystyle\max_{t\in \R}|\{x\in G: g(x)=t\}|$. For any sets $\A\subset G$ and $\B,\C\subset \R^*$, we have 
$$|f(\A,\B)||\B+\C|\gg \min\Bigg\{\dfrac{q^r|\B|}m,\dfrac{|\A||\B|^2|\C|}{m^2q^{2r-1}}\Bigg\}.$$
\end{theorem}
Here, and throughout, $X\gg Y$ means that there exists a positive constant $c$ such that $X\geq cY$, and $"\ll"$ is defined in a similar way. We also say that $X\sim Y$ if $X\gg Y$ and $Y\gg X$. These above results are generalizations of earlier results due to Hegyv\'{a}ri and Hennecart \cite{heg} in the prime field setting. A prime field version for small sets can also be found in the work on Mojarrad and Pham \cite{pham2}.

For any $\A\subset \R^*$, it follows from Theorem \ref{ham} with $g(x)=x$ and $h(x)\equiv 0$ that

\[ |\A\cdot \A||\A+\A|\gg  \min \left\lbrace q^{r} |\A|, ~\frac{|\A|^4}{q^{2r-1}}\right\rbrace.\]
Thus,
\[\max\left\lbrace |\A+\A|, |\A\cdot\A|\right\rbrace \gg  \min \left\lbrace q^{\frac{r}{2}} |\A|^{\frac{1}{2}}, ~\frac{|\A|^2}{q^{\frac{2r-1}{2}}}\right\rbrace.\]

For $\A\subset\R$, we define $\A^2:=\{x^2\colon x\in \A\}$. In a recent work, Yazici \cite{ay} studied another sum-product type estimate. Namely, she proved the following theorem. 
\bigskip
\begin{theorem}[Yazici, \cite{ay}]\label{thmz}
Let $\R$ be a finite valuation ring of $q^r$ and $\A$ be a subset of $\R$. If $|\A+\A||\A|^2\gg q^{3r-1}$, then 
\[\max\left\lbrace |\A+\A|, |\A^2+\A^2|\right\rbrace \gg  q^{\frac{r}{4}} |\A|^{\frac{3}{4}}.\]
\end{theorem}
It follows from Theorem \ref{thmz} that either the size of $\A+\A$ or $\A^2+\A^2$ is large when $|\A|$ is big enough. Note that this sum-product type estimate was first studied by Solymosi in \cite{solymosi}.  In this paper, we provide two improvements of this result. Our first result is stated as follows. 
\bigskip
\begin{theorem}\label{thm0}
Let $\R$ be a finite valuation ring of order $q^r$ with $q$ a power of an odd prime number. For $\A\subset \R$ with $|\A|\ge 2q^{r-1}$, then 
\[\max\left\lbrace |\A+\A|, |\A^2+\A^2|\right\rbrace \gg  \min \left\lbrace q^{\frac{r}{2}} |\A|^{\frac{1}{2}}, ~\frac{|\A|^2}{q^{\frac{2r-1}{2}}}\right\rbrace.\]
\end{theorem}
One can check that Theorem \ref{thm0} improves Theorem \ref{thmz} when $|\A|\gg q^{r-\frac{1}{3}}$. Indeed, it is clear that 
\[\min \left\lbrace q^{\frac{r}{2}} |\A|^{\frac{1}{2}}, ~\frac{|\A|^2}{q^{\frac{2r-1}{2}}}\right\rbrace=q^{\frac{r}{2}}|\A|^{\frac{1}{2}}.\]
Thus, under the condition $|\A|\gg q^{r-\frac{1}{3}}$, the conclusion of Theorem \ref{thm0} is stronger than that of Theorem \ref{thmz}. We also need to compare the conditions of these two theorems. Obviously, the assumption $|\A|\gg q^{r-\frac{1}{3}}$ implies the condition $|\A+\A||\A|^2\gg q^{3r-1}$ by using the fact that $|\A+\A|\ge |\A|$. 
When $|\A+\A||\A|^2\gg q^{3r-1}$ and $|\A|\ll q^{r-\frac{1}{3}}$, we have another improvement of Theorem \ref{thmz}, which is a consequence of the following theorem. 
\\
\begin{theorem}\label{theorem1}
Let $\R$ be a finite valuation ring of order $q^r$ with $q$ a power of an odd prime number. For $\A\subset \R$ with $|\A|\ge 2q^{r-1}$ and $|\A+\A||\A|^2\gg q^{3r-1}$, then
$$\max \big\{|\A+\A|, |\A^2+\A^2|\big\}\gg q^{r/3}|\A|^{2/3}.$$
\end{theorem}
It is convenient to give a brief comparison between these theorems. It follows from Theorem \ref{thm0} that if $|\A|\ll q^{r-\frac{1}{3}}$, then $\max\left\lbrace |\A+\A|, |\A^2+\A^2|\right\rbrace \gg q^{\frac{1-2r}{2}}|\A|^2$, which is better than the bound $q^{\frac{r}{4}}|\A|^{\frac{3}{4}}$ of Theorem \ref{thmz} whenever $|\A|\gg q^{r-\frac{2}{5}}$, and weaker than the threshold $q^{r/3}|\A|^{2/3}$ of Theorem \ref{theorem1} whenever $|\A|\ll q^{r-\frac{3}{8}}$. Thus, in the range $q^{r-\frac{3}{8}}\ll |\A|\ll q^{r-\frac{1}{3}}$, the bound 
$\frac{|\A|^2}{q^{\frac{2r-1}{2}}}$ is the best. If $|\A+\A||\A|^2\gg q^{3r-1}$ and $2q^{r-1}\le |\A|\ll q^{r-\frac 38}$, then the lower bound of Theorem \ref{thmz} is the strongest. In other words, we can summarize the bounds in the following corollary.
\\
\begin{corollary}\label{hiep}
Let $\R$ be a finite valuation ring of order $q^r$ with $q$ a power of an odd prime number, and $\A$ be a set in $\R$. 
\begin{enumerate}
\item If $|\A|\gg q^{r-\frac{1}{3}}$, then 
\[\max\left\lbrace |\A+\A|, |\A^2+\A^2|\right\rbrace \gg  q^{\frac{r}{2}}|\A|^{\frac{1}{2}}.\]
\item If $q^{r-\frac{3}{8}}\ll |\A|\ll q^{r-\frac{1}{3}}$, then
\[\max\left\lbrace |\A+\A|, |\A^2+\A^2|\right\rbrace \gg  \frac{|\A|^2}{q^{\frac{2r-1}{2}}}.\]
\item If $|\A+\A||\A|^2\gg q^{3r-1}$ and $2q^{r-1}\le|\A|\ll q^{r-\frac{3}{8}}$, then
\[\max\left\lbrace |\A+\A|, |\A^2+\A^2|\right\rbrace \gg  q^{r/3}|\A|^{2/3}.\]
\end{enumerate}
\end{corollary}
Furthermore, we have a remark on the last statement of Corollary \ref{hiep}. If $|\A+\A||\A|^2\gg q^{3r-1}$ and $|\A|\ll q^{r-\frac{3}{8}}$, then $|\A+\A|\gg |\A|q^{\frac{1}{8}}$. This leads to $\A+\A$ is an expanding set. However, the lower bound $q^{r/3}|\A|^{2/3}$ for $\max\left\lbrace |\A+\A|, |\A^2+\A^2|\right\rbrace$ is stronger whenever $|\A|\ll q^{r-\frac{3}{8}}$. We refer the interested readers to \cite{bou, tv} for related sum-product results in the finite ring setting.

The main difference between our method and that of Yazici is that she used the Pl\"{u}nnecke-Ruzsa inequality. Instead of which, to prove Theorem \ref{thmz} and Theorem \ref{thm0}, we will use spectral graph theory techniques and some ideas from the work of Pham, Vinh, and De Zeeuw \cite[Theorem 1.3]{t}. Especially, our method can be easily extended to the case of higher dimensions. In particular, we obtain the following main results which are extensions of Theorem \ref{thm0} and Theorem \ref{theorem1}, respectively. (Here, for a positive integer $n$, we use notation $n\A^2$ for the set consisting of all elements of the form $a_1+a_2+\cdots+a_n$ with $a_1,a_2,\ldots, a_n$ in $\A^2$.)
\\
\begin{theorem}\label{thm1}
Let $\R$ be a finite valuation ring of order $q^r$ with $q$ a power of an odd prime number. For any $\A\subset \R$ with $|\A|\ge 2q^{r-1}$ and for any integer $n>1$, we have
\[\max\left\lbrace |n\A^2|, ~|\A+\A|\right\rbrace\gg \min \left\lbrace q^{\frac{r}{n}}|\A|^{\frac{n-1}{n}}, ~\frac{|\A|^{\frac{3n-2}{n}}}{q^{\frac{(n-1)(2r-1)}{n}}}\right\rbrace.\]
\end{theorem}
\begin{theorem}\label{thm2}
Let $\R$ be a finite valuation ring of order $q^r$ with $q$ a power of an odd prime number. For any $\A\subset \R$ with $|\A|\ge 2q^{r-1}$ and for any integer $n>1$, we have
\[\max \{|\A+\A|, |n\A^2|\} \gg q^{\frac{r}{2n-1}}|\A|^{\frac{2n-2}{2n-1}},\]
whenever $|\A+\A|^{n-1}|\A|^n\gg q^{r+(n-1)(2r-1)}$.
\end{theorem}
In the rest of the paper, we are going to give proofs of Theorem \ref{thm1} and Theorem \ref{thm2}.
\section{The definition of finite valuation rings}
We start this section by recalling the definition of finite valuation rings from \cite{n}.
\bigskip

A commutative ring with identity is called a {\it finite valuation ring} if it is finite, local and principal.

Let $\mathcal R$ be a finite valuation ring. We have that $\R$ contains a unique maximal ideal, denoted by $(z)$ for some non-unit \textit{uniformizer} element $z$ in $\R$. Notice that the uniformizer element is defined up to a unit in $\R$. We also denote by $\R^*, \R^0$ the set of units, non-units in $\R$, respectively.

Since $(z)$ is the maximal ideal, we have $\R/(z)$ is a field, which is denoted by $F$. We denote the size of $F$ by $q$. Let $r$ be the nilpotency degree of $z$, i.e. the smallest positive integer $r$ satisfying $z^r=0$. It is known that $q$ is a power of a prime number. In this paper, we assume that $q$ is odd, so $2$ is a unit in $\R$, i.e $2\in \R^*$.

Over finite valuation rings $\R$,  one has a natural valuation 
$$f: \R\rightarrow\{0,1,\ldots, r\},$$
defined by $f(0)=r$ and for $x\neq 0$, $f(x)=k$ if $x\in (z^k)\setminus (z^{k+1}).$ This means that, $f(x)=k$ if and only if $x=uz^k$ for some unit $u$ in $\R$. One can check that for each $k$, the group group $(z^k)/(z^{k+1})$ is a one-dimensional linear space over the residue field $F=\R/(z)$, thus its size is $q$. Hence, $|(z^k)|=q^{r-k}$ for $k=0,1,\ldots,r$. In other words, we have $|(z)|=q^{r-1}, |\R|=q^r$ and $|\R^*|=|\R|-|(z)|=q^r-q^{r-1}$. We refer the readers to \cite{n} for more details. There are several examples of finite valuation rings, for instance, finite fields $\mathbb{F}_q$   with $q$ is a prime power, finite cyclic rings $\mathbb{Z}/p^r\mathbb{Z}$ with $p$ is a prime, and $\mathcal O/(p^r)$, where $\mathcal O$ is the ring of integers in a number field and $p\in\mathcal O$ is a prime.
\section{Techniques from spectral graph theory}
Let $G=(A\cup B, E)$  be a bipartite graph. If all vertices in each part have the same degree, we say that $G$ is \textit{biregular}. If $G$ is biregular, we denote the common degree of each vertex in $A$ by $\deg (A)$, and the degree of each vertex in $B$ by $\deg(B)$. Let $M$ be the adjacency matrix of $G$. Assume that $\lambda_1, \ldots, \lambda_n$ are eigenvalues of $M$ with $|\lambda_1|\geq |\lambda_2|\geq\cdots \ge|\lambda_n|$. Since $G$ is a bipartite graph, we have $\lambda_2=-\lambda_1$. To prove our main theorems, we will make use of the following version of the Expander mixing lemma for bipartite graphs. We refer the readers to \cite{e} for a proof.
\bigskip

\begin{lemma}\label{lm1}
Suppose that $G=(A\cup B, E)$ is biregular with $\deg(A)=a$ and $\deg(B)=b$.  For subsets $X\subset A$ and $Y\subset B$, the number of edges between $X$ and $Y$, denoted by $e(X, Y)$, satisfies
$$\Big|e(X,Y)-\dfrac a {|B|}|X||Y|\Big|\leq \lambda_3\sqrt{|X||Y|},$$
where $\lambda_3$ is the third eigenvalue of $G$. We note that $a/|B|=b/|A|$.
\end{lemma} 

For any $\x\in\R^d\setminus {(\R^0)^d}$,
 denote by $[\x]$ the equivalence class of $\x$ in $\R^d\setminus{(\R^0)^d}$, 
where $\x, \y\in \R^d\setminus {(\R^0)^d}$ are in the same class iff $\x=t\y$ for some unit $t\in \R^*$. Let $\E_{q,d}(\R)=(A\cup B, E)$ be  the Erd\H{o}s-R\'enyi bipartite graph with $A$ and $B$ being the sets of equivalence classes in  $\R^d\setminus{(\R^0)^d}$. There is an edge between two vertices $[\x]$ and $[\y]$ iff $\x\cdot \y=0$.  The spectrum of this graph over finite valuation rings was studied by Nica \cite{n} by using exponential sums. We summary in the following theorem.
\bigskip
\begin{theorem}[{\bf Nica}, \cite{n}]\label{thr} The cardinality of each vertex part of 
$\E_{q,d}(\R)$ is $q^{(d-1)(r-1)}(q^d-1)/(q-1)$, and 
$\deg(A)=\deg(B)=q^{(d-2)(r-1)}(q^{d-1}-1)/(q-1)$. 
The third eigenvalue of $\E_{q,d}(\R)$ is at most $\sqrt{q^{(d-2)(2r-1)}}.$
\end{theorem}

\section{Proof of Theorem \ref{thm1}}
Since $|\A|\geq 2q^{r-1}$, it follows that 
$$|\A\cap \R^*|\geq |\A|-|\R^0|=|\A|-q^{r-1}\geq \dfrac {|\A|}2.$$
Thus we may assume that $\A$ is a subset of $\R^*$. We now prove that the size of the set $\A^2$ is at least $\gg |\A|$. Indeed, suppose $x^2=y^2$ with $x,y\in A$, then $(x-y)(x+y)=0$. There are posibilities for pairs $(x,y)$ as follows: $x=y$, or $x=-y$, or $x-y\in (z)$ and $x+y\in (z)$. If the last case holds, then we can write $x-y=u_1z^k_1$ and $x+y=u_2z^{k_2}$ with $u_1, u_2\in\R^*$ and some positive integers $k_1, k_2$. This leads to $2x=u_1z^{k_1}+u_2z^{k_2}\in (z)$, which gives a contradiction since both $2 $ and $x$ are in $\R^*$. In other words, either $x=y$ or $x=-y$ and therefore $|\A^2|\gg |\A|.$

Define $\D=n\A^2$. Consider the following equation
\[ x+(b_1-c_1)^2+\cdots +(b_{n-1}-c_{n-1})^2=t, \]
where $x\in \A^2, b_i\in \A+\A, c_i\in \A, 1\le i\le n-1, t\in D$. Let $N$ be the number of solutions of this equation. We first see that $N\ge |\A|^{2n-1}$. Let $U$ and $V$ be two vertex sets of the Erd\H{o}s-R\'{e}nyi graph $\mathcal{E}_{q, n+1}(\R)$ defined by
\[U:=\Bigg\{ \bigg(-2b_1, \ldots, -2b_{n-1}, \sum_{i=1}^{n-1}b_i^2+x, 1\bigg)\colon b_i\in \A+\A, 1\le i\le n-1, x\in \A^2\Bigg\},\]
and 
\[V:=\Bigg\{\bigg (c_1, \ldots, c_{n-1}, 1, \sum_{i=1}^{n-1}c_i^2-t \bigg)\colon c_i\in \A, 1\le i\le n-1, t\in \D\Bigg\}.\]
We have $|U|\sim |\A+\A|^{n-1}|\A|$ and $|V|=|\A|^{n-1}|\D|$. It is not hard to check that $N$ is bounded by the number of edges between $U$ and $V$ in the graph $\mathcal{E}_{q, n+1}(\R)$. Thus, one can apply Lemma \ref{lm1} and Theorem \ref{thr} to get 
\[N\ll \frac{|U||V|}{q^r}+q^{\frac{(n-1)(2r-1)}{2}}\sqrt{|U||V|}=\frac{|\A+\A|^{n-1}|\A|^n|n\A^2|}{q^r}+q^{\frac{(n-1)(2r-1)}{2}}\sqrt{|\A+\A|^{n-1}|\A|^n|n\A^2|}.\]
Using the fact that $N\ge |\A|^{2n-1}$, we obtain 
\[\max\left\lbrace |n\A^2|, ~|\A+\A|\right\rbrace\gg \min \left\lbrace q^{\frac{r}{n}}|\A|^{\frac{n-1}{n}}, ~\frac{|\A|^{\frac{3n-2}{n}}}{q^{\frac{(n-1)(2r-1)}{n}}}\right\rbrace.\]
This completes the proof of the theorem.
\hfill $\square$

\section{Proof of Theorem \ref{thm2}}
Since $|\A|\geq 2q^{r-1}$, as in the proof of Theorem \ref{thm1}, we may assume that $\A$ is a subset of $\R^*$.  In this proof, we will follow the idea of \cite[Theorem 1.3]{t}.

As in the proof of Theorem \ref{thm1}, we define $D=n\A^2$. 
Consider the following equation
\[ x+(b_1-c_1)^2+\cdots +(b_{n-1}-c_{n-1})^2=t, \]
where $x\in \A^2, b_i\in \A+\A, c_i\in \A, 1\le i\le n-1, t\in D$. Let $N$ be the number of solutions of this equation. We see that $N\ge |\A|^{2n-1}$.

By the Cauchy-Schwarz inequaltiy, one has 
\[N^2\le |\D|\cdot E,\]
where $E$ is the number of tuples $(x, b_1, \ldots, b_{n-1}, c_1, \ldots, c_{n-1}, y, d_1, \ldots, d_{n-1}, e_1, \ldots, e_{n-1})$ satisfying 
\[x+(b_1-c_1)^2+\cdots +(b_{n-1}-c_{n-1})^2=y+(d_1-e_1)^2+\cdots +(d_{n-1}-e_{n-1})^2.\]
Let $U$ and $V$ be two vertex sets in the Erd\H{o}s-R\'{e}nyi graph $\mathcal{E}_{q, 2n}(\R)$ defined by
\[U:=\left\lbrace \bigg(-2b_1, \ldots, -2b_{n-1}, 2d_1, \ldots, 2d_{n-1}, 1, \sum_{i=1}^{n-1}b_i^2-\sum_{i=1}^{n-1}d_i^2+x\bigg)\colon b_i\in \A+\A, d_i\in \A, x\in \A^2\right\rbrace,\]
and 
\[V:=\left\lbrace \bigg(c_1, \ldots, c_{n-1}, e_1, \ldots, e_{n-1}, \sum_{i=1}^{n-1}c_i^2-\sum_{i=1}^{n-1}e_i^2-y, 1\bigg)\colon c_i\in \A, e_i\in \A+\A, y\in \A^2\right\rbrace.\]
We have $|U|=|V|\sim |\A+\A|^{n-1}|\A|^n$. We also have $E$ is bounded by the number of edges between $U$ and $V$ in the graph $\mathcal{E}_{q, 2n}(\R)$. Therefore, it follows from Lemma \ref{lm1} and Theorem \ref{thr} that 
\[E\ll \frac{|\A+\A|^{2n-2}|\A|^{2n}}{q^r}+q^{(n-1)(2r-1)} |\A+\A|^{n-1}|\A|^n.\]
Using the facts that $N\ge |\A|^{2n-1}$ and $N^2\le |\D|\cdot E$, we derive 
\[\max \{|\A+\A|, |n\A^2|\} \gg q^{\frac{r}{2n-1}}|\A|^{\frac{2n-2}{2n-1}},\]
whenever $|\A+\A|^{n-1}|\A|^n\gg q^{r+(n-1)(2r-1)}$. In particular,
\[\max \{|\A+\A|, |n\A^2|\} \gg \min \left\lbrace q^{\frac{r}{2n-1}}|\A|^{\frac{2n-2}{2n-1}}, \frac{|\A|^{\frac{3n-2}{n}}}{q^{\frac{(n-1)(2r-1)}{n}}}\right\rbrace,\]
which ends the proof of the theorem. 
$\hfill \square$

\noindent{\sc Duc Hiep Pham}\\
University of Education\\Vietnam National University, Hanoi\\ 144 Xuan Thuy, Cau Giay, Hanoi\\ Vietnam\\
Email: {\sf phamduchiep@vnu.edu.vn}

\end{document}